\documentclass[10pt,french]{article}

\usepackage{amssymb,amsmath,amsthm, amscd}
\usepackage[francais]{babel}
\usepackage[latin1]{inputenc}
\usepackage[all]{xy}
\usepackage{amssymb}
\usepackage{pstricks}
\usepackage{amsmath}
\usepackage{array}
\usepackage{indentfirst}
\usepackage{hyperref}
\usepackage{graphicx}
\usepackage[latin1]{inputenc}
\usepackage{amssymb}
\usepackage{pstricks}
\usepackage{amsmath}
\usepackage{array}
\usepackage{indentfirst}
\usepackage{hyperref}
\usepackage{graphicx}

\newcommand{\F}{{\mathbb F}_2}

\newcommand{\R}{\mathbb{R}}

\newcommand{\map}{\mathrm{map}}
\newcommand{\U}{\mathcal U}
\newcommand{\Z}{\mathbb Z}
\newcommand\ra{\rightarrow}

\newcommand\C{{\mathcal C}}

\newcommand{\A}{\mathcal A}

\newcommand{\bqn}{\begin{eqnarray*}}
\newcommand{\eqn}{\end{eqnarray*}}

\def \B{\mathrm{B}}
\def \C{\mathrm{C}}
\def \F{\mathrm{F}}

\def \N{\mathrm{N}}

\def \R{\mathrm{R}}
\def \Z{\mathbb{Z}}

\def \colim{\mathrm{colim}}

\def \hom{\mathrm{hom}}

\def \I{\mathcal{I}}
\def \Im{\mathrm{Im}}

\def \fp{\mathbb{F}_p}

\def\U{\mathcal{U}}

\def\hatY{\widehat{Y}}

\def \underline#1{\underset{\widetilde{~}}{#1}}

\def\mappt{\mathrm{map}_*}

\title{Applications depuis $K(\Z/p,2)$ et une conjecture de N. Kuhn}

\author{G. Gaudens, L. Schwartz}

\begin{document}

\maketitle

\begin{abstract}
On d\'emontre une conjecture due \`a N. Kuhn concernant la
cohomologie singuli\`ere \`a coefficients mod $p$ des espaces, comme module instable sur l'alg\`ebre de Steenrod. Notre
d\'emonstration de ce r\'esultat, d\'ej\`a connu en caract\'eristique
$2$, fait appel \`a une m\'ethode  nouvelle, qui
fonctionne en toute caracteristique. De cette mani\`ere on
r\'etablit un r\'esultat de \cite{S98} dont la preuve est incompl\`ete dans le cas d'un nombre premier impair.
\end{abstract}

\renewcommand{\abstract}{ \begin{center} {\bf Abstract.} \end{center}}

\begin{abstract}
We settle a  conjecture due to N. Kuhn about the mod $p$ cohomology of spaces considered as unstable modules over the Steenrod algebra. This result is already known to hold  in characteristic $2$. The method presented here is essentially new and works for all characteristics. In doing so we fix a gap in \cite{S98} concerning the odd prime case.
\end{abstract}

\tableofcontents

\oddsidemargin=-0.5cm
\evensidemargin=-0.5cm
\topmargin=-1.5cm
\parindent=0pt
\textheight=25 cm
\textwidth=16 cm

\newtheorem{thm}{Th\'eor\`eme}[section]
\newtheorem{lem}[thm]{Lemme}
\newtheorem{cor}[thm]{Corollaire}
\newtheorem{prop}[thm]{Proposition}
\newtheorem{defi}[thm]{D\'efinition}
\newtheorem{rmq}[thm]{Remarque}
\def\exe{{\bf Exemples}}
\def\B{{\mathcal B}}
\theoremstyle{definition} \def\Ker{{\rm Ker}}
\def\Im{{\rm Im}}
\newtheorem{ex}[thm]{Exemple}

\theoremstyle{remark}

\newtheorem{rem}[thm]{Remarque}
\newtheorem{preuve}[thm]{Preuve}
\def\resp{(resp. $\C$) }

\numberwithin{equation}{section}
\def\R{\mathbb R}
\def\N{\mathbb N}
\def\F{{\mathbb F}_p}
\def\U{\mathcal U}
\def\Z{\mathbb Z}
\def\rg{{\rm rang}}

\def\A{{\mathcal A}_p}
\def\NK{{\mathcal N \mathcal K}}

\def\An{{\bf An.}}

\def\I{{\bf Inv.}}

\begin{section}{Introduction, la conjecture  de Kuhn.}

Soit $p$ un nombre premier. La cohomologie $H^* (X; \F)$ que l'on notera $H^*X$  d'un  espace  $X$ est naturellement un objet de la cat\'egorie  $\U$ des modules instables sur l'alg\`ebre de Steenrod \cite{S94}.
Cette cat\'egorie est munie, comme toute cat\'egorie ab\'elienne,  d'une  filtration naturelle, dite de Krull,  par des sous-cat\'egories \'epaisses stables par colimites :
$$
\U _0 \subset \U _1 \subset \U _2  \subset \ldots \subset \U
$$
La sous-cat\'egorie $\U _0$ est exactement la sous-cat\'egorie des
modules localement finis \cite{S94} obtenue comme colimite de
modules finis. Autrement dit cette sous-cat\'egorie est obtenue
\`a partir des objets simples de $\U$, les modules $\Sigma^n\F$,
par  extensions et colimites. Les termes suivants sont d\'efinis
it\'erativement comme suit: $\U _n$ est `l'image inverse' dans
$\U$ de  la sous-cat\'egorie ab\'elienne $\mathcal C$ du quotient
$\U/\U_{n-1}$ obtenue \`a partir des objets simples de
$\U/\U_{n-1}$ par  extensions et colimites. La cat\'egorie $\U$
\emph{n'est pas} la r\'eunion des $\U_n$.

La filtration de Krull est caract\'eris\'ee par :

\begin{thm}\cite{S94}
\label{theorem:filtration:krull} {\it Soit $\bar T$ le foncteur de
lannes r\'eduit, adjoint \`a gauche du produit tensoriel par $\tilde
H^*(B\Z/p)$. Un module instable $M$ est dans $\U_n$ si et
seulement si $\bar T^{n+1} (M)=\{0\}$.}
\end{thm}

Un bref rappel de la th\'eorie de Lannes appara\^{\i}t en section \ref{cohomologie:des:espaces:d:applications:pointees}.
Tout module instable $M$ est \'egalement muni de sa  filtration
nilpotente $\{M_s\}_{s\geq 0}$ qui est decroissante, s\'epar\'ee
et naturelle, et dont les sous-quotients sont les suspensions
it\'er\'ees de modules instables r\'eduits, plus pr\'ecis\'ement:
$$
M_s/M_{s+1}= \Sigma ^s R_s M
$$
o\`u $R_sM$ est r\'eduit (voir \cite{S94}).
L'objet du pr\'esent article est de d\'emontrer le r\'esultat suivant conjectur\'e par N. Kuhn.

\begin{thm}
\label{theoreme:principal} {\it Soit $X$ un espace. Si la cohomologie $H^* X$ est dans un cran fini de la filtration de Krull, alors $H^* X$ est  localement finie.
En d'autres termes, si $H^* X$ est dans $\U_n$ pour un certain $n$, alors $H^* X$ est dans $\U_0$.}
\end{thm}

Si cet \'enonc\'e ne fait appara\^itre que la filtration de Krull,
on verra que la d\'emonstration en repose sur les rapports entre
filtration de Krull et filtration nilpotente.

Afin de donner du relief \`a ce r\'esultat, observons que la
cat\'egorie $\U$ poss\`ede une infinit\'e d'objets qui sont dans
un cran fini de la filtration de Krull et qui ne sont pas
localement finis, par exemple les g\'en\'erateurs projectifs
$F(n)$ et leurs suspensions \cite{S94}; et qu'\`a l'inverse,    la
cohomologie modulo $p$ d'un espace d'Eilenberg-Mac Lane $K(\Z/p ,
n)$ pour $n>0$, ainsi que la cohomologie modulo $p$ d'un groupe
fini (d'ordre divisible par $p$), sont des modules instables qui
n'appartiennent \`a $\U _n$ pour aucun entier $n$. Ce r\'esultat
\'etait d\'ej\`a connu pour $p=2$ \cite{S01, DG03}. Il corrige une
erreur qui    concerne le cas d'un nombre premier impair dans
l'article \cite{S98}   du second auteur. On reviendra plus loin
sur ce point. L'approche pr\'esent\'ee ici est  nouvelle et unifie
le cas $p=2$ avec le cas o\`u $p$ est  impair. Ce r\'esultat
implique  le corollaire suivant (\'egalement conjectur\'e par
Kuhn) :

\begin{cor}
\label{corollaire:theoreme:principal} {\it Soit $X$ un espace.  Si la cohomologie $H^* X$ a une famille finie de g\'en\'erateurs comme module instable sur l'alg\`ebre de Steenrod, alors $H^* X$ est de dimension finie comme $\fp$-espace vectoriel.}
\end{cor}

Observons que ce r\'esultat est essentiellement instable: la cohomologie modulo $p$ d'un spectre d'Eilenberg-Mac Lane $H\Z/p$ est monog\`ene comme module  sur l'alg\`ebre de Steenrod, mais n'est pas de dimension finie comme $\fp$-espace vectoriel.

On d\'erive le Corollaire \ref{corollaire:theoreme:principal} du
Th\'eor\`eme \ref{theoreme:principal}  en observant que la
cat\'egorie $\U$ poss\`ede un syst\`eme de g\'en\'erateurs
projectifs, et que chacun d'entre eux est dans un cran fini de la
filtration de Krull \cite{S94}. Il s'ensuit qu'un module instable
finiment engendr\'e est dans un cran fini de la filtration de
Krull, car celle-ci est une filtration croissante par des
sous-cat\'egories \'epaisses. Suivant  le Th\'eor\`eme
\ref{theoreme:principal}, si la cohomologie $H^* X$ a un nombre
fini de g\'en\'erateurs comme module instable, alors  elle est
localement finie. Mais un module finiment engendr\'e et localement
fini doit \^etre de dimension finie comme $\fp$-espace vectoriel.

\vskip 5mm

L'article est organis\'e comme suit. Dans la fin de l'introduction
on explique o\`u se situe l'erreur dans \cite{S98}. Dans la
premi\`ere partie on rappelle bri\'evement \emph{la r\'eduction de
Kuhn du probl\`eme}. On se donne un module instable $M$
satisfaisant aux hypoth\`eses requises, et \'etant en particulier
la cohomologie d'un espace $X$; la cohomologie de
$X$ est de plus supppos\'ee $s$-nilpotente ($s>0$) exactement. Puis on
construit une application alg\'ebrique, r\'ealisable
topologiquement par une application $\varphi$,  de $H^*(\Omega
\Sigma X)$ vers $\Sigma^s \tilde H^*B\Z/p$. On montre ensuite que
l'application $\varphi$ ne peut s'\'etendre \`a $\Sigma^{s-1}
K(\Z/p,2)$ le long de l'application naturelle $\Sigma^{s} B\Z/p
\to B\Z/p$, car l'action des op\'erations de Steenrod ne serait
pas respect\'ee.

Dans les deux parties finales on montre que les groupes
d'obstructions \`a l'extension sont nuls, ce qui m\`ene \`a la
contradiction. La derni\`ere construction est  r\'eminiscente de
celle de H. Miller \`a la fin de son article sur la conjecture de
Sullivan aux Annals \cite{Mil84}.

\vskip 5mm

La preuve donn\'ee pour le cas o\`u $p$ impair  dans \cite[Section
3]{S98}  contient une erreur. L'approche de \cite{K08} ne traite
que le cas $p=2$ et pr\'esenterait une difficult\'e analogue \`a
celle de \cite{S98} pour $p$ impair d'apr\`es M. Stelzer. L'objet
du pr\'esent article est donc \'egalement de montrer comment
\'eviter le probl\`eme  rencontr\'e dans \cite{S98}.

Il   est  expliqu\'e dans cet article pourquoi  la m\'ethode utilis\'ee pour le cas $p=2$ ne peut l'\^etre
pour $p$ impair sans modification.  L'analogue (impair) des classes consid\'er\'ees pour traiter du cas $p=2$  se trouvent  sur la colonne $-p$ de la suite spectrale d'Eilenberg-Moore et  peuvent  donc supporter une
diff\'erentielle $d_{p-1}$ non triviale, pour $p=2$ il n'y a \'evidement pas de diff\'erentielle! En fait on a  la formule (voir \cite{Sm})
$$
d_{p-1}(x \otimes \ldots \otimes x)= \lambda \beta P^k(x)
$$
avec $k=\frac{\vert x \vert-1}{ 2}$. Si l'homomorphisme de
Bockstein est trivial dans le module consid\'er\'e la m\'ethode
fonctionne comme pour $p=2$ et redonne le r\'esultat de \cite{K95}
sans faire appel au th\'eor\`eme sur l'invariant de Hopf  ou \`a
d'autres r\'esultats profonds de th\'eorie de l'homotopie. En fait
ceci  \'etend m\^eme le th\'eor\`eme de Kuhn de \cite{K95} qui
montre que pour qu'un module (ou une alg\`ebre instable)  du  type
consid\'er\'e soit (\'eventuellement) r\'ealisable il faut que des
homomorphismes de Bockstein soient non triviaux entre deux
\'etages de la filtration nilpotente du module consid\'er\'e, la
m\'ethode ci dessus montre que ce doit \^etre le cas entre les
deux premiers \'etages non triviaux.

Par contre la m\'ethode propos\'ee dans \cite{S98} pour
`\'eliminer' le Bockstein :   prendre le plan
projectif sur l'espace des lacets est inexacte, et  \'echoue sans
argument suppl\'ementaire. Il reste qu'en examinant de plus pr\`es
des cas particuliers on voit rapidement appara\^itre des
contradictions qui forcent la nullit\'e de l'homomorphisme de
Bockstein et permettent donc de se ramener au cas pr\'ec\'edent.

A titre d'exemple  consid\'erons le cas o\`u l'alg\`ebre instable
initialement consid\'er\'e a un quotient non-tivial dans
${\mathcal N}il_1$ et que tous les produits sont nuls. Cette
derni\`ere condition est facile \`a assurer en remplacant l'espace
suppos\'e exister par un quotient. Soit $x$ la classe
consid\'er\'ee plus haut, de degr\'e $2k+1$ et $y=P^kx$, $k=p^h$
pour un certain $h$. Si l'homomorphisme de Bockstein est non
trivial sur $x$ il l'est sur $y$ car le module consid\'er\'e est
de la forme $A \otimes \Phi^{h+1} F(1)$ (voir \cite{S98}).
Rappelons que $F(1)$ est le module instable librement engendr\'e
en degr\'e $1$, qui s'identifie au module des primitifs de  la
coalg\`ebre $H^* B\Z/p$.

L'\'el\'ement $y ^{\otimes p-1}  \otimes x$ existe sur la colonne
$-p$ du terme $E_2$ de la suite spectrale et ne peut \^etre tu\'e
par une diff\'erentielle $d_i$, $i \leq p-1$ pour des raisons de
degr\'es, quitte encore une fois \`a prendre un quotient de
l'espace initial. Toutes les diff\'erentielles sont nulles sur cet
\'el\'ement (en particulier il n'y a pas de place pour des
produits de Massey non triviaux, voir \cite{Cleary85}, chapitres 7 et 8). On a, comme plus haut,
$$
d_{p-1}(y \otimes \ldots \otimes y)= \lambda \beta P^{pk}(y)
$$
et  donc :
$$
0=P^kd_{p-1}( y \otimes \ldots \otimes y  \otimes x)=d_{p-1}P^k( y \otimes \ldots \otimes y  \otimes x)= \lambda \beta P^{pk}(y) \not = 0
$$
On a donc une contradiction si $\lambda \not =0$.

Cependant, pour \'etablir le r\'esultat en toute g\'en\'eralit\'e,
il convient de proc\'eder autrement, comme il va d\'ecrit
ci-dessous.

Le corollaire  est d\'emontr\'e pour $p=2$   via la
non-r\'ealisabilit\'e comme cohomologie d'un espace topologique de
certains modules instables finis dans \cite{S98} et \cite{K08}).
Cette approche ne fonctionne en l'\'etat que pour $p=2$, il n'est
pas imm\'ediatement clair comment appliquer les techniques
ci-dessous aux complexes finis bien quelques faits sugg\`erent que
des adaptations puissent \^etre trouv\'ees.

\vskip 4mm

Dans ce qui suit on fixe un nombre premier $p$ et on suppose que les espaces sont $p$-complets. Ceci n'entra\^ine pas de perte de g\'en\'eralit\'e. On supposera aussi que l'homologie modulo $p$ des espaces consid\'er\'es est
de dimension finie en chaque degr\'e et que cette propri\'et\'e persiste si on applique le foncteur $T$ de Lannes. On montrera dans une derni\`ere section comment s'affranchir de cette condition en utilisant la th\'eorie des espaces profinis suivant \cite{DG03}.

\end{section}

\section{Construction d'applications et th\'eorie d'obstruction}

On revient dans cette section sur la strat\'egie de
d\'emonstration du th\'eor\`eme. On \'etablit le r\'esultat par
l'absurde: on suppose qu'il existe un espace dont la cohomologie
n'est pas localement finie mais appartient \'a un cran fini de la
filtration  de Krull $\U_n$ avec $n>0$. Dans \cite{K95} Kuhn
montre que l'on peut se ramener au cas o\`u $n=1$. On suppose donc
donn\'e un espace dont la cohomologie est dans $\U_1$, mais pas  dans
$\U_0$. On montre alors qu'il existe une application de cet espace
(ou d'un espace d\'eduit de cet espace) vers une suspension d'un
espace d'Eilenberg-MacLane $K(\Z/p,2)$ telle que l'application
induite en cohomologie ne commute pas aux op\'erations de
Steenrod.

On ne rappelle pas ici en d\'etail comment se ramener \`a supposer
que le module est dans $\U_1$, ceci est fait dans \cite{K95} et
repris dans les articles suivants. Simplement pour montrer la
filiation de la construction ci-dessous avec celle de Kuhn
observons que son outil central est l'application du th\'eor\`eme
de Lannes sur la cohomologie des espaces fonctionnels de source
$B\Z/p$ \`a la cofibre de $X \ra \map(B\Z/p,X)$, dont la
cohomologie sous des hypoth\`eses favorables, est donn\'ee par
$\bar T(H^*X)$.

On  commence par montrer comment construire des applications
r\'ealisant certaines fl\`eches alg\'ebriques.

\subsection{Application alg\'ebrique}

Soit $X_s$ un espace $s$-connexe dont la cohomologie r\'eduite est
exactement $s$-nilpotente avec $s>0$ et non localement finie. On
suppose de plus que $R_s H^*X_s$ est dans $\U_1$. Un module
r\'eduit  dans $\U_1$ est la somme d'un module concentr\'e en
degr\'e $0$ et d'un sous-module non trivial d'une somme directe 
$\oplus _{\lambda \in \Lambda} F(1)$  \cite{S98}.

En s\'electionnant  un $\lambda_0$ dans l'ensemble d'indices
$\Lambda$ de la somme on a une application alg\'ebrique non
triviale
$$
\varphi^*_s: H^* X_s \to \Sigma ^s R_s H^* X_s \to \Sigma ^s F(1)  \subset \Sigma^s  H^* B\Z/p~.
$$

Dans la sous-section qui suit, on montre que l'application
$\varphi_s^* $ est induite par une application d'espaces :

\begin{lem}Il existe une application $\varphi_s :  \Sigma ^s
B\Z/p \to X_s $ induisant l'application $\varphi_s^*$ en
cohomologie modulo $p$.
\end{lem}

\subsection{Construction de l'application $\varphi_s$}
\label{construction}

Nous donnons deux  solutions pour l'existence d'une telle
application, voici la premi\`ere. 

On suppose donc  $R_s H^* X_s $ est dans $\U_1$, avec $H^* X_s$ exactement $s$-nilpotent. Posons $X_t
=\Omega^{s-t}X_s$. D'apr\`es \cite[Th\'eor\`emes 1.5 et 1.6]{S01}
l'espace $ \Omega^{s-1} X_s$ poss\`ede une cohomologie r\'eduite
exactement $1$-nilpotente, et l'application naturelle $R_s H^* X_s
\ra  R_1 H^* X_1 $ a un noyau et un conoyau nilpotent. Ce
`F-isomorphisme' selon la terminologie de \cite{S01} est induit
par l'\'evaluation
$$
\Sigma^{s-1}\Omega^{s-1}
X_s = \Sigma^{s-1}  X_1 \to X_s ~.
$$

De $\varphi_s^*$ on d\'eduit alors $\varphi_1^*$

$$
\varphi^*_1: H^* X_1 \to \Sigma ^1 R_1 H^* X_1 \to \Sigma  F(1)  \subset \Sigma  H^* B\Z/p~.
$$

Ceci  ram\`ene au cas o\`u on veut construire
$$
\varphi_1 : \Sigma
B\Z/p \to X_1= \Omega ^{s-1} X_s~,
$$
car si on y parvient, l'application adjointe
$$
\varphi_s : \Sigma^s B\Z/p \to X_s
$$
jouit des propri\'et\'es d\'esir\'ees.
Rappelons maintenant que $H^*\Sigma B\Z/p $ est $\U$-injectif. Il
s'ensuit, que le morphisme de Hurewicz
$$
\map (\Sigma B\Z/p , X) \to \hom _{\mathcal{K}} (H^* X ,  H^* \Sigma B\Z/p)
$$
est surjectif \cite{LS89, GL87}. Ici $\mathcal{K}$ d\'esigne la
cat\'egorie des alg\`ebres instables \cite{S94}. Ceci assure
l'existence d'une application $\varphi_1$, une fois observ\'e que
$\varphi_1^*$ est un morphisme dans $\mathcal{K}$.

\vskip 2mm

Que l'application alg\'ebrique $\varphi_s^*$ soit r\'ealisable
topologiquement r\'esulte \'egalement de l'argument  ci dessous.

On suppose que $X_s$ est un $H$-espace et que l'on peut calculer
la cohomologie de $\map _*(B\Z/p,X_s)$ comme il est indiqu\'e en
section 3. En particulier la cohomologie de cet espace est
exactement $(s-1)$-connexe. On peut choisir une classe dans le
$s$-i\`eme groupe d'homologie et la repr\'esenter par une
application de $S^s$ dans $\map_*(B\Z/p, X_s)$. On consid\`ere
l'application adjointe $S^s \times B\Z/p \ra X_s$. Elle est
triviale sur $S^s \times *$, car c'est l'adjointe d'une
application \`a valeurs dans l'espace des applications point\'ees.
Elle est homotopiquement triviale sur $S^s \times *$ car on peut supposer sans restriction que  $X_s$ est au moins $2s$-connexe.
Par le th\'eor\`eme de Lannes sur les applications de source $B\Z/p$ \cite{La92}, puisque $\tilde H^*X_s$ est nilpotente, elle est homotopiquement triviale  sur $* \times B\Z/p$.

On r\'ecup\`ere donc une application $\Sigma^s  B\Z/p= S^s \wedge B\Z/p \ra X_s$ qui a l'action prescrite en cohomologie.

\subsection{Impossibilit\'e d'\'etendre $\varphi_s$ et cons\'equences}

\begin{lem}L'application $\varphi_s:  \Sigma^s B\Z/p
\to X_s$ ne s'\'etend pas en une application $\Sigma ^{s-1}K(\Z/p,
2) \to X_s$ le long de l'application $\Sigma^{s-1}\psi: \Sigma^s
B\Z/p \to \Sigma^{s-1} K(\Z/p, 2)$.
\end{lem}

En effet, l'existence d'une telle extension  permettrait de factoriser  le morphisme de modules instables non trivial   $\varphi_s^* :\tilde H^*X_s \ra \Sigma \tilde H^*B\Z/p$  au travers de
$ \tilde H^* \Sigma ^{s-1} K(\Z/p, 2)
$
ce qui est  impossible, car il n'y a pas d'applications non triviales d'une suspension $s$-i\`eme vers une suspension $(s-1)$-i\`eme d'un module r\'eduit.
\vskip 2mm

L'espace $K(\Z/p,2)$ se reconstruit \`a partir $B\Z/p$ \`a l'aide de la construction de Milnor. On a une filtration $* =C_0 \subset C_1= \Sigma B\Z/p \subset C_2 \subset \ldots \subset \cup_n C_n=K(\Z/p,2)$, et un diagramme o\`u les fl\`eches de la ligne sup\'erieure sont homotopiquement triviales (et $B$
d\'esigne $B\Z/p$) :
\begin{equation*}
\begin{CD}
\cdots@>{}>>  B^{*n} @>{}>>B^{*n+1}@>{}>>\cdots\\
@VVV @VVV @VVV @ VVV \\
\cdots@>{}>>  C_{n} @>{}>>C_{n+1} @> >> \cdots\\
\end{CD}
\end{equation*}
d'o\`u on d\'eduit des cofibrations \`a homotopie pr\`es :
$$
\Sigma^{n-1}B^{\wedge n}  \ra  C_n \ra C_{n+1} \ra \Sigma^n B^{\wedge n}
$$

Puisque l'application consid\'er\'ee plus haut ne peut s'\'etendre \`a $\Sigma^{s-1} K(\Z/p, 2)$ on a n\'ecessairement :

\begin{lem} Il existe un entier $n$ tel que le groupe
$[\Sigma ^{n+s-2}  (B\Z/p)^{\wedge n}, X_s]$ est non trivial.
\end{lem}

On a $[\Sigma ^{n+s-2}  (B\Z/p)^{\wedge n}, X_s]= \pi_{n+s-2}
\map_* ( B\Z/p^{\wedge n},  X_s)$.

\begin{lem}
\label{lemme:non:trivial}  Il existe un entier  $n$ tel que $H^{*}
\map_* ( B\Z/p^{\wedge n}, X_s)$  n'est pas $(n+s-2)$-connexe.
\end{lem}

\begin{section}{Cohomologie des espaces d'applications point\'ees}
\label{cohomologie:des:espaces:d:applications:pointees}

Soit $H$ la cohomologie modulo $p$ de l'espace classifiant $B\Z/p$, et soit $\bar H$ sa cohomologie r\'eduite.  Les endofoncteurs de la categorie $\U$
$$
M\mapsto M \otimes H; ~~~~ M \mapsto M \otimes \bar H
$$
poss\`edent des adjoints \`a gauche respectifs $T$ et $\bar T$.
Les foncteurs $T$ et $\bar T$ sont tous deux exacts. De plus  $T$
commute aux produits tensoriels dans le sens o\`u
l'application naturelle $$T (M_1)\otimes T (M_2)\to T
(M_1\otimes M_2)$$
donn\'ee formellement par les propri\'et\'es
d'adjonction, est un isomorphisme pour tous $M_1$ et $M_2$. Le
scindement naturel $H\cong \F \oplus  \bar H$ dans la cat\'egorie
$\U$ conduit \`a un scindement naturel
$$
TM \cong \bar T M \oplus M .
$$

Supposons l'espace $X$ $p$-complet, $1$-connexe, et que $TH^*X$ soit de dimension finie en chaque degr\'e. Alors, \cite{La92}, l'application  naturelle :
$$
T H^*X \to H^* \map (B\Z/p,X)
$$
adjointe de l'application induite en cohomologie par
l'\'evaluation $B\Z/p \times \map (B\Z/p, X) \to X$, est un
isomorphisme d'alg\`ebres  instables. De plus le facteur $\bar T
H^*  \subset T H^* X$ identifie $\bar T H^* X$ avec la cohomologie
r\'eduite de $\Delta X$, la cofibre  homotopique de l'application
naturelle $X\subset \map (B\Z/p, X)$ donn\'ee par les applications
constantes.

Soit $X$ un H-espace, typiquement un espace de lacets.
On a alors un scindement
$$
\map (B\Z/p, X)\cong  \map _* (B\Z/p, X) \times X.
$$
Suivant \cite{CCS07} on en  d\'eduit:

 \begin{lem}{\it $Q H^* \map _*
(B\Z/p^{\wedge n}, X )= \bar{T}^n QH^* X$.}
\end{lem}

La connexit\'e de l'id\'eal d'augmentation d'une alg\`ebre connexe
co\"{\i}ncide avec celle de ses ind\'ecomposables. En
particulier, on d\'eduit du Lemme  \ref{lemme:non:trivial}:

\begin{prop}\label{nc} Supposons que l'application consid\'er\'ee ci dessus ne s'\'etende pas \`a tout $\Sigma^{s-1} K(\Z/p, 2)$ et  que $X_s$ soit un H-espace.  Alors il
existe un $n$ tel que $QH^{*} \map_* ( B\Z/p ^{\wedge n+1}, X_s)= \bar T
^ {n} Q H^* X_s $ ne soit  pas $(n+s-2)$-connexe.
\end{prop}

On va montrer que dans le cas  consid\'er\'e ceci n'a pas lieu, ce qui ach\`evera la d\'emonstration.

\begin{rmq}
En g\'en\'eral on ne sait pas calculer la cohomologie de l'espace
des applications point\'ees.  Il se trouve que sous les
hypoth\`eses sur le module $A$ est de la forme $\Phi^k F(1)
\otimes F$, $F$ fini, et que l'on consid\`ere (s'il en existe) un
espace $p$-complet $X$ tel que $\tilde H^*X \cong A$ on peut
calculer $H^* \map _*(B\Z/p,X)$ comme module instable. On commence
par supposer que les conditions d'application du th\'eor\`eme de
Lannes sont satisfaites (et que l'entier $k$ est grand
relativement \`a $F$).  Notons $Y$ la cofibre de l'application $X
\ra \map (B\Z/p,X)$. On montre alors que $\map (B\Z/p,X)$ est, \`a
$p$-compl\'etion pr\`es, homotopiquement \'equivalent au bouquet
$X \vee Y$.  Il r\'esulte alors, par exemple, de  \cite{Strom03}
que la suspension de l'espace d'applications point\'ees est
homotopiquement \'equivalent  \`a $\Sigma Y \vee (\Sigma Y \wedge
\Omega X)$. On n'\'echappe pas donc \`a calculer la cohomologie de
$\Omega X$ (de laquelle on tirera la contradiction). Il vaut donc
mieux consid\'erer la d'embl\'ee comme il est fait ci-dessous le cas
de l'espace des fonctions dans un espace de lacets.
\end{rmq}

\end{section}

\begin{section}{Conjecture de Kuhn}

\subsection{R\'eduction du probl\`eme}

Supposons donn\'e un espace  $Z$ tel que $H^* Z$ n'est pas
localement finie, mais appartienne \`a un cran fini $\U _n$ de la
filtration de Krull. Le module $R_0 \tilde H^* Z$ est
n\'ecessairement trivial, car dans le cas contraire, l'existence
d'une classe non-r\'eduite contredit le fait que $H^* Z$ est de
filtration de Krull finie (ce point est amplement document\'e dans
les articles pr\'ec\'edents sur le sujet, ou simplement dans
\cite{S94}; cela suit egalement de \cite{GS05}. En utilisant les r\'eductions donn\'ees dans
\cite{K95} (voir aussi \cite{DG03}), on obtient qu'il existe alors
un espace connexe $Z'$ dont la cohomologie r\'eduite est
exactement $s$-nilpotente avec $s>0$ et non localement finie. De
plus on peut supposer la cohomologie de  $X$ est $2s$-connexe et $R_s H^*X$ est
dans $\U_1$. On pose alors $X_s:= \Omega \Sigma X$.

\subsection{Connexit\'e de $QH^* X_s$ et contradiction}

La structure de $H^* X_s$ comme module instable,  en fonction de
$H^* X$ est tr\`es classique et donn\'ee par le th\'eor\`eme de
Bott-Samelson (voir \cite{GW}) :

\begin{thm}  Soit  $K:= H^* X_s$, c'est une alg\`ebre de Hopf. Elle est isomorphe  \`a l'alg\`ebre tensorielle $\mathbb{T}H^* X$
\begin{itemize}
\item comme $\fp$-espace vectoriel,
\item comme module instable,
\item  comme coalg\`ebre  avec le coproduit de d\'econcat\'enation,
\end{itemize}
La filtration par la longueur des tenseurs sur $K\cong \mathbb{T}H^* X$ est multiplicative, et induit le  produit de battage (`shuffle product´) sur le gradu\'e.
\end{thm}
 
Rappelons que $\tilde H^*X_s$ est suppos\'ee $s$-nilpotente. L'alg\`ebre de cohomologie est le produit tensoriel d'alg\`ebres ext\'erieures et d'alg\`ebres polynomiales tronqu\'ees \`a hauteur $p$. Pour une telle alg\`ebre instable, on va montrer:

\begin{lem}\label{alg}{\it Pour tout $n\geq 1$,
$\bar T^{n+1} QK$ est $((n+1)s-1)$-connexe}.
\end{lem}

Pour $s>0$, l'in\'egalit\'e $(n+1)s-1\geq n+s-1$ produit une
contradiction entre la Proposition \ref{nc} et le Lemme  \ref{alg}, ce qui \'etablit le
th\'eor\`eme. Il reste donc \`a d\'emontrer \ref{alg} .

On  consid\`ere la sous-alg\`ebre  $K_h$ de $K $
engendr\'ee par le sous-module $M_h$ des  tenseurs de longueur non-nulle 
inf\'erieure ou \'egale  \`a $h$. On a donc  par construction pour $h\geq 1$ un \'epimorphisme
$$
M_h \to QK_h
$$
Comme $M_h$ est dans $\U_h$, il en est de m\^eme pour $QK_h$,
de sorte que $\bar T^{h+1}QK_h=0$.

On a ${\mathrm{colim}_h }  \, Q K_h = QK$, et comme $\bar T$ commute aux colimites, on a
$$
\bar T ^{n+1} Q K \cong \bar T ^{n+1}Q (\colim_h K_h)\cong \colim_h
(\bar T^{n+1} Q K_h) \cong  \bar T^{n+1} \colim_{h>n} QK_h
$$

On a un diagramme commutatif
$$
\xymatrix{ M_h \ar[r]\ar[d]& M_{h+1} \ar[r]\ar[d]& M_{h+1}/M_h \ar[d]\ar[r]&\{0\}\\
QK_h \ar[r] & Q K_{h+1} \ar[r]& Q (K_{h+1}/(\bar K_h.{K}_{h+1})) \ar[r]&\{0\}}
$$
dont les fl\`eches verticales sont des \'epimorphismes, et dont les
suites horizontales sont exactes, et ce pour tout $h\geq 1$.
En appliquant le foncteur exact $\bar T^{n+1}$ \`a suite exacte
horizontale inf\'erieure du diagramme pour $h=n$, on voit que
$$
\bar T^{n+1} Q K_{ n+1} \cong \bar T^{n+1} Q
(K_{n+1}/(\bar{K}_n.K_{n+1}))~.
$$
$ M_{n+1}/M_ n$ s'identifie comme module instable avec $(\tilde
H^* X )^{\otimes n+1}$- qui est $(n+1)s$-nilpotent, donc
$((n+1)s-1)$-connexe. Il suit que $\bar T^{n+1} Q K_{ n+1}$ est
$((n+1)s-1)$-connexe. Ensuite pour $h\geq n+1$, en utilisant la
suite exacte
$$
\bar  T^{n+1} QK_h \to  \bar T^{n+1} Q K_{h+1} \to  \bar  T^{n+1} Q
(K_{h}/(K_h.\bar{K}_{h+1})) \ra \{0\}
$$
et le fait que $\bar T^{n+1}Q (K_{h+1}/(K_h.K_{h+1}))$ est un
quotient de $\bar T^{n+1}(H^* X )^{\otimes h+1}$- qui est $((h+1)s-1)$,
connexe on voit que $\bar T^{n+1}QK_{h+1}$ est
$((n+1)s-1)$-connexe si $\bar T^{n+1}QK_h$ l'est.

Par r\'ecurrence sur $h\geq n+1$ et passage \`a la colimite, on en d\'eduit que $ \bar T^{n+1} QK_h$ est
$((n+1)s-1)$-connexe, comme annonc\'e.
\end{section}

\begin{section}{De l'utilisation des espaces profinis}

Nous donnons dans cette derni\`ere section, comme annonc\'e en introduction,  les \'el\'ements de la th\'eorie des espaces profinis qui nous permettent de travailler sans hypoth\`ese de finitude, et qui n'apparaissent pas d\'ej\`a  dans \cite{MO, DG03}- auxquels nous renvoyons pour plus de  d\'etails.

Rappelons qu'un ensemble profini est un espace topologique muni d'une topologie qui le rend compact et compl\'etement discontinu. La cohomologie modulo $p$  $H^* X$d'un espace profini $X$, d\'efinie dans \cite{MO}, est naturellement munie d'une structure d'alg\`ebre instable sur l'alg\`ebre de Steenrod.
Morel a montr\'e que les isomorphismes en cohomologie continue modulo $p$ et les monomorphismes sont les \'equivalences faibles et les cofibrations d'une cat\'egorie de mod\`eles simpliciale sur la cat\'egorie des espaces profinis, pour laquelle il existe  un remplacement fibrant fonctoriel $Y\to \hatY$ (voir \cite{DG03}). Le foncteur d'oubli de la topologie des espaces profinis dans les espaces ordinaires poss\`ede un adjoint \`a gauche  $X\to \widehat{X}$ appel\'e \emph{compl\'etion profinie}. L'application naturelle  $H^* Y\to H^*\hatY$ est un \emph{isomorphisme} d'alg\`ebres instables, donc de modules instables. En utilisant les remarques qui suivent sur les espaces fonctionnels et la connexit\'e pour les espaces profinis, le lecteur se convaincra qu'on montre le r\'esultat suivant, en rempla\c{c}ant  \emph{espace} par \emph{espace profini} dans le corps de l'article.

\begin{thm} {\it Soit $X$ un espace profini. Si la cohomologie $H^* X$ est dans un cran fini de la filtration de Krull, alors $H^* X$ est  localement finie.
En d'autres termes, si $H^* X$ est dans $\U_n$ pour un certain $n$, alors $H^* X$ est dans $\U_0$.}
\end{thm}

Avec l'observation que la cohomologie d'un espace ordinaire est canoniquement isomorphe \`a celle de sa comp\'etion profinie, on obtient les r\'esultats
\'enonc\'es en introduction  \emph{sans hypoth\`ese de finitude}.

Pour tout espace profini $Y$ et tout ensemble fini simplicial, on peut former un espace profini fonctionnel $\map (X, Y)$, limite inverse du syst\`eme $\map (X, \hatY (-))$, o\`u $\hatY (-)$ d\'esigne le remplacement fibrant fonctoriel de \cite{MO}, voir \'egalement \cite{DG03}. On n'aura uniquement \`a consid\'erer le cas o\`u $X= \Sigma ^s B\Z/p^{\wedge n}$, qui est  bien un ensemble fini simplicial. Pour  $X$ et $Y$ des espaces profinis point\'es, on a la notion d'un espace profini fonctionnel point\'e  $\map_*(X, Y)= \lim \map_* (X, \hatY (-))$. Dans le cas o\`u $X= S^1$, on obtient $\Omega Y := \map_* (S^1, Y)$. On peut d\'efinir la suspension d'un espace profini point\'e $X= \lim X(-)$ comme $\Sigma X := \lim \Sigma X(-)$; la formule d'adjonction usuelle $\hom (\Sigma Z, Y) \cong\hom ( Z, \Omega Y) $ est alors valide, o\`u $\hom$ d\'esigne l'ensemble des morphismes entre deux espaces profinis.
Le foncteur $\mappt (-, Y)$ est une fibration. En particulier $\map_*(- , Y)$ transforme les suites de Puppe de cofibrations  en suite de Puppe de fibrations. Pour tout ensemble profini point\'e $Y$, on peut introduire la notion d'ensembles/groupes d'homotopie profinis $\pi_n (Y):= \lim \pi_n \hatY (-)$. Ces derniers sont munis d'une topologie  profinie naturelle. Une fibration entres espaces point\'es induit une suite exacte longue sur les groupes d'homotopie, tout comme dans le cas usuel.
Soit $[X, Y]$ l'ensemble des classes d'homotopie simpliciale d'applications continues point\'ees  $X \to
\hatY$. On a une formule $[X, Y ]:=  \pi_0\map_* (X, Y)$ et il s'ensuit que:
$$
[\Sigma^n X , Y]\cong \pi_0 \map_*(\Sigma_n X, Y)\cong \pi_0
\map_*(X, \Omega ^n Y)\cong \pi_0 \Omega^n \map_* (X, Y)\cong \pi_n \map _*
(X, Y)
$$
Nous notons enfin la version profinie du th\'eor\`eme de Hurewicz:
\begin{thm}Soit $n\geq 0$ et soit $X$ un espace profini connexe point\'e.
Alors $\pi_i X $ est trivial pour tout $i\leq n$ si et seulement si $\tilde{H} ^i
X$ est trivial pour tout $i\leq n$. Enfin, si $\pi_l X$ est le groupe d'homotopie de $X$ non trivial de plus bas degr\'e, alors $H^l(X)$ s'identifie au groupe des homomorphismes continus  $\pi_l X \to \Z/p$.

\end{thm}

\end{section}

\bibliographystyle{amsalpha}

\providecommand{\bysame}{\leavevmode\hbox to3em{\hrulefill}Thinspace}
\providecommand{\MR}{\relax\ifhmode\unskip\space\fi MR }
\providecommand{\MRhref}[2]{%
  \href{http://www.ams.org/mathscinet-getitem?mr=#1}{#2}
}
\providecommand{\href}[2]{#2}

\end{document}